\documentclass[11pt]{amsart}
\usepackage{mathrsfs}
\usepackage{amsfonts}
\usepackage{amssymb}
\usepackage{marvosym, manfnt,  yhmath,  stmaryrd}
\usepackage{graphicx}
\usepackage[]{epsfig}
\usepackage{amssymb}
\usepackage[]{pstricks}
\newpsobject{malla}{psgrid}{subgriddiv=1,griddots=10,gridlabels=6pt}
\newcommand{\lan}{\langle}
\newcommand{\ran}{\rangle}
\newtheorem{theorem}{Theorem}[section]
\newtheorem{proposition}[theorem]{Proposition}

\newtheorem{lemma}[theorem]{Lemma}
\newtheorem{remark}[theorem]{Remark}

\newcommand{\cF}{{\mathcal F}}

\newcommand{\cX}{{\mathcal X}}

\newcommand{\R}{{\mathbb R}}

\title[Critical Schr\"odinger-Debye System]{Local and Global Well-Posedness for the Critical Schr\"odinger-Debye System}

\author{Ad\'an  J. Corcho}

\author{Filipe Oliveira}

\author{Jorge Drumond Silva}

\address{\textbf{Ad\'an J. Corcho}
\newline
Instituto de Matem\'atica.
\newline
Universidade Federal do Rio de Janeiro-UFRJ.
\newline
Ilha do Fund\~{a}o, 21945-970. Rio de Janeiro-RJ, Brazil.
\newline
Rio de Janeiro-RJ, Brazil.} \email{adan@im.ufrj.br}

\address{\textbf{Filipe Oliveira}
\newline
Centro de Matem\'atica e Aplica\c{c}\~oes, FCT-UNL
\newline
Monte da Caparica, Portugal}
\email{fso@fct.unl.pt }

\address{\textbf{Jorge Drumond Silva}
\newline
Centro de An\'alise Matem\'atica, Geometria e Sistemas Din\^amicos,
\newline
Departamento de Matem\'atica,
\newline
Instituto Superior T\'ecnico,
\newline
Av. Rovisco Pais, 1049-001 Lisboa, Portugal.}
\email{jsilva@math.ist.utl.pt }

\thanks{A. J. Corcho was supported by CAPES and CNPq (Edital Universal-482129/2009-3), Brazil\\
J. Drumond Silva was partially supported by Funda\c c\~ao para a Ci\^encia e Tecnologia (FCT/Portugal) through
the program PCI 2010/FEDER\\Filipe Oliveira was partially supported by FCT/Portugal through Financiamento Base 2008-ISFL-
1-297.}



\begin{document}
\maketitle

\setcounter{page}{1}

\begin{quote}
{\normalfont\fontsize{8}{10}
\selectfont {\bfseries Abstract.}
We establish local well-posedness results for the Initial
Value Problem associated to the Schr\"odinger-Debye system in dimensions $N=2, 3$
for data in $H^s\times H^{\ell}$, with $s$ and $\ell$ satisfying  $\max \{0, s-1\} \le \ell \le \min\{2s, s+1\}$. In particular, these include the energy space $H^1\times L^2$. Our results improve the previous ones obtained in \cite{Bidegaray1}, \cite{Bidegaray2} and \cite{Corcho-Linares}. Moreover, in the critical case ($N=2$) and for initial data in $H^1\times L^2$, we prove that solutions exist for all times, thus providing a negative answer to the open problem mentioned in \cite{Fibich-Papanicolau} concerning the formation of singularities for these solutions.
\par}
\end{quote}

\section{\textbf{Introduction}}
We consider the Initial Value Problem (IVP) for the Schr\"odinger-Debye system
\begin{equation}\label{SD}
\begin{cases}
iu_t+\frac{1}{2}\Delta u=uv, & t\ge0,\; x\in\mathbb{R}^{N}\; (N=1,2,3),\\
\mu v_t +v = \lambda|u|^{2}, & \mu>0,\; \lambda=\pm1,\\
u(x,0)=u_{0}(x),\quad v(x,0)=v_{0}(x), &
\end{cases}
\end{equation}
where $u=u(x,t)$ is a complex-valued function, $v=v(x,t)$ is a real-valued
function and $\Delta$ is the Laplacian operator in the spacial variable.
This model describes the propagation of an electromagnetic wave through a nonresonant medium whose material response
time is relevant.
See Newell and Moloney \cite{Newell} for a more complete
discussion of this model.

In the absence of delay ($\mu =0$), the system (\ref{SD}) reduces to the cubic Nonlinear Schr\"odinger equation
(\textbf{NLS})
\begin{equation}\label{NLS}
iu_t+\tfrac{1}{2}\Delta u=\lambda u\left\vert u\right\vert ^2,
\end{equation}
which is focusing or defocusing for $\lambda=-1$ and $\lambda=1$, respectively. Similarly,
the sign of the parameter $\lambda$  provides an analogous  classification of (\ref{SD}).

For sufficiently regular data, the \emph{mass} of the solution $u$ of the system (\ref{SD}) is invariant. More
precisely,
\begin{equation}
\int_{\R^N}|u(x,t)|^{2}dx=\int_{\R^N}|u_{0}%
(x)|^{2}dx.\label{Conservation Law}%
\end{equation}
Other conservation laws for this system are not known, but
the following pseudo-Hamiltonian structure holds:
\begin{equation}\label{Energy-1}
\frac{d}{dt}E(t)=2\lambda\mu\int_{\R^N}(v_t)^{2}dx,
\end{equation}
where
\begin{equation}\label{Energy-2}
E(t)=\int_{\R^N}\bigl\{|\nabla u|^{2}+\lambda|u|^{4}-\lambda\mu^{2}(v_t)^2\bigl\}dx
=\int_{\R^N}\bigl\{|\nabla u|^{2}+2v|u|^2-\lambda v^2\bigl\}dx.
\end{equation}

The system (\ref{SD}) can be decoupled  by solving the second equation with respect to $v$,
\begin{equation}\label{eq. de v}
v(t)= e^{-t/\mu}v_0(x)+ \tfrac{\lambda}{\mu}\int\limits_0^t\,
e^{-(t-t')/\mu}|u(t')|^2\,dt',
\end{equation}
to obtain the integro-differential equation
 \begin{equation}\label{SD-IDF}
 \begin{cases}
 iu_t+\tfrac{1}{2}\Delta u = e^{-t/\mu}uv_0(x)
+ \tfrac{\lambda}{\mu}u\int\limits_0^t\,
 e^{-(t-t')/\mu}|u(t')|^2dt',&  x\in {\R}^N,\; t\ge 0, \\
 u(x,0) = u_0(x).
 \end{cases}
 \end{equation}

 The rest of this introduction is organized as follows: in section 1.1 we review the previous existing results regarding the local and global theory for (\ref{SD}). In section 1.2, we describe our new results in dimensions $N=2,3$.
\subsection{Overview of former results in dimensions $\boldsymbol {N=1, 2, 3}$}

$\,$\vspace{0.3cm}

We begin with a review of the
local and global theory for the Cauchy problem (\ref{SD}) with initial data $(u_{0},v_{0})$ in Sobolev spaces
$H^{s}(\R^N)\times H^{\ell}(\R^N)$,\, $N=1,2,3$.

Bid\'egaray (\cite{Bidegaray1} and  \cite{Bidegaray2})
established the following local results:
\begin{theorem}[\emph{Bid\'egaray, 2000}]\label{LWP-Bidegaray}
Let $N=1,2,3$ and $(u_0,v_0)\in H^s(\R^N)\times H^s(\R^N)$.  The IVP (\ref{SD-IDF}) has a unique solution
\begin{enumerate}
\item [(a)] $u \in L^{\infty}\left([0,T]; H^s(\R^N)\right)$\; if $s>N/2$,
\smallskip
\item [(b)] $u \in L^{\infty} \left ([0,T]; H^1(\R^N)\right)$\; if $s=1$,
\smallskip
\item [(c)] $u \in C\left([0,T]; L^2(\R^N)\right) \bigcap L^{8/N}\left ([0,T]; L^{4}(\R^N)\right)$\; if $s=0$,
\end{enumerate}
where $T=T(\|u_0\|_{H^s}, \|v_0\|_{H^s})>0$. Moreover, the solution $u$ depends continuously on the initial data
$(u_0,v_0)$.
\end{theorem}
\noindent
These results were obtained by a fixed-point procedure applied to the Duhamel formulation
for the integro-differential equation (\ref{SD-IDF}), using the Strichartz estimates for the unitary Schr\"odinger group
\begin{equation}\label{Schrodinger-Group}
S(t)f(x)= \left (\text{e}^{-it|\xi|^2/2}\widehat{f}(\xi)\right)^{\lor}(x).
\end{equation}

Following the same approach,  Corcho and Linares (\cite{Corcho-Linares}) improved
the results stated in Theorem \ref{LWP-Bidegaray} in the one-dimensional case. More precisely, they
established the fo\-llowing assertions:
\begin{theorem}[\emph{Corcho-Linares, 2004}]\label{LWP-Corcho-Linares}
Let $0 < s \le 1$,  $q \in {[2,\infty]}$,\; $2/r=1/2-1/q$\, and\, $(u_0,v_0) \in H^s(\R)\times H^{\ell}(\R)$.
The IVP (\ref{SD-IDF}) has a unique solution
\begin{enumerate}
\item [(a)] $u \in \cX^{s,r, q}_{T}:=C([0,T];\;H^s(\R))\bigcap L^r([0,T]; L^q(\R))$\, if\, $0< s < 1/2$ and $\ell =s$,
\item [(b)] $u \in \cX^{s,r, q}_{T}$ and $u_x\in L^{\infty}\left (\R; L^2([0,T])\right)$\, if\,  $s=1/2$ and $0\le \ell \le 1/2$,
\item [(c)] $u \in \cX^{s,r, q}_{T}$ and $u_x\in L^{\infty}\left (\R; L^2([0,T])\right)$\,  if\, $1/2<s\le 1$ and $s-1/2 < \ell \le s$,
\end{enumerate}
where $T=T(\mu, \|u_0\|_{H^s}, \|v_0\|_{H^{\ell}})>0$. Moreover, the map\; $(u_0,v_0) \longmapsto u(t)$ is locally Lipschitz and  $v\in C\bigl([0,T];\;H^{\ell}(\R)\bigl)$.
\end{theorem}
\noindent The new ingredients used in the proof of Theorem \ref{LWP-Corcho-Linares} are commutator estimates for fractional Sobolev spaces and
the smoothing effect for the Schr\"odinger group
$$\|D^{1/2}S(t)u_0\|_{L^{\infty}_xL^2_T}\le C\|u_0\|_{L^2},$$
deduced by Kenig, Ponce and Vega (see \cite{KPV0-a,KPV0-b}).
Furthermore, the authors also showed that, although the fixed point procedure is performed only on the function $u$,
equation (\ref{eq. de v}) can be used to obtain the persistence property of the solution $v$ in $H^s(\R^N)$ in the
cases described in Theorem \ref{LWP-Bidegaray}.

Concerning global existence,
it was also proved in \cite{Corcho-Linares} that
the local-in-time results for the solution $u$
of the integro-differential equation (\ref{SD-IDF}), given in Theorem \ref{LWP-Bidegaray} (c) and
Theorem \ref{LWP-Corcho-Linares} (b) and (c), can be extended to all positive times.
However, the method used does not provide control of
the evolution in time of the $H^s$-norm of the corresponding solution $v$. Indeed, contrarily to the NLS equation,
(\ref{SD}) does not possess a Hamiltonian structure, hence the extension to any positive times
of the local-in-time solutions $(u,v)$
is not straighforward.

Recently, however,
Corcho and Matheus  (see \cite{Corcho-Matheus}) studied the case $N=1$ in the framework of Bourgain spaces
and obtained the following local and global well-posedness results for the system:
\begin{theorem}[\emph{Corcho-Matheus, 2009}]\label{LWP-Corcho-Matheus}
For any $(u_0,v_0)\in H^s(\R) \times H^{\ell}(\R)$, where
\begin{equation}\label{local-theorem-continuous-a}
|s|-1/2\le \ell < \min\{s+1/2,\; 2s+1/2\}\quad and \quad  s> -1/4,
\end{equation}
there exists a time $T=T(\|u_0\|_{H^{s}}, \|v_0\|_{H^{\ell}})>0$ and
a unique solution $(u(t),v(t))$ of the initial value problem
(\ref{SD}) in the time interval $[0,T]$, satisfying
$$(u, v)\in C\left([0,T]; H^s(\R) \times H^{\ell}(\R)\right).$$
Moreover, the map $(u_0,v_0) \longmapsto (u(t),v(t))$ is locally
Lipschitz. In addition, in the case $\ell =s$ with $-3/14 < s \le 0$, the
local solutions can be extended to any time interval $[0,T]$.
\end{theorem}
\noindent The global results in Theorem \ref{LWP-Corcho-Matheus} are based on a good control of the
$L^2$-norm of the solution $v$, which provides
global well-posedness in $L^2\times L^2$. Global well-posedness below $L^2$ regularity is then obtained
via the \emph{I-method} introduced by Colliander, Keel, Staffilani, Takaoka and Tao in \cite{Imethod}.

Regarding the formation of singularities in the critical case ($N=2$), Fibich and Papanicolaou
(\cite{Fibich-Papanicolau}) studied this system in the focusing case using the lens transformations, but did not
derive any result as to the blow-up of the solutions. On the other hand, from a numerical point of view,
Besse and Bid\'egaray (\cite{Besse-Bidegaray}) used two different methods suggesting that
blow-up occurs for initial data $u_0(x,y)=\text{e}^{-(x^2+y^2)}$ and $v_0=\lambda|u_0|^2$.
However, prior to the present paper, the blow-up problem remained open.

\subsection{Main results in dimensions $\boldsymbol{N=2, 3}$}

$\,$\vspace{0.3cm}

In this paper we give a negative answer to the question of existence of blow-up solutions for initial data in $H^1(\R^2)\times L^2(\R^2)$ (see Theorem \ref{global-theorem}). Note that this result is not in contradiction with the numerical simulations in \cite{Besse-Bidegaray}. Indeed, in the latter, the suggested blow-up occurs for the norm $\|u(\cdot,t)\|_{L^{\infty}}$ which, in two dimensions, is not controlled by $\|u(\cdot,t)\|_{H^1}$. Also, contrarily to the NLS case, we prove that the blow-up occurs neither in the defocusing nor in the focusing case. This is due to the delay induced by the term $\mu v_t$ in the left-hand-side of the second equation of (\ref{SD}), which prevents the solution from concentrating critically.
As expected, this behavior does not depend on the size of $\mu$, as long as this parameter stays positive.
This was already remarked in \cite{Besse-Bidegaray}: for if $(u,v)$ is a solution to (\ref{SD}) for a value of $\mu >0$, then
$\Bigl(\tilde{u}(x,t),\tilde{v}(x,t)\Bigl)=\left(\mu^{1/2}u(\mu^{1/2}x, \mu t),\,\mu v(\mu^{1/2}x, \mu t)\right)$
yields a solution to (\ref{SD}) for $\mu =1$.

\begin{figure}[h]\label{Figure-LWP}
\centering \psset{unit=1cm}
\begin{pspicture}(-2,-2)(4,4)%
\malla%
\psline{->}(-2,0)(4,0)%
\psline{->}(0,-2)(0,4)%
\pspolygon[fillstyle=solid,fillcolor=cyan,linewidth=0.5pt]
(0,0)(1,0)(4,3)(4,4)(3,4)(1,2)(0,0)%
\rput(0,0){$\bullet$}%
\rput(1,0){$\bullet$}%
\rput(2,2){$\mathcal{W}$}%
\rput(3.7,-0.25){$s$}%
\rput(0.25,3.8){$\ell$}%
\rput(1.5,3.5){\small{$\ell=s+1$}}%
\rput(3.2,1.2){\small{$\ell=s-1$}}%
\rput(0.2,1.5){\small{$\ell=2s$}}%
\end{pspicture}
\caption{{\small\it{The region $\mathcal{W}$,
bounded by the lines $\ell=0$, $\ell=s-1$, $\ell =2s$ and $\ell=s+1$, corresponds to the set of indices $(s,\ell)$ of our local well-posedness results for system (\ref{SD}).}}}
\end{figure}
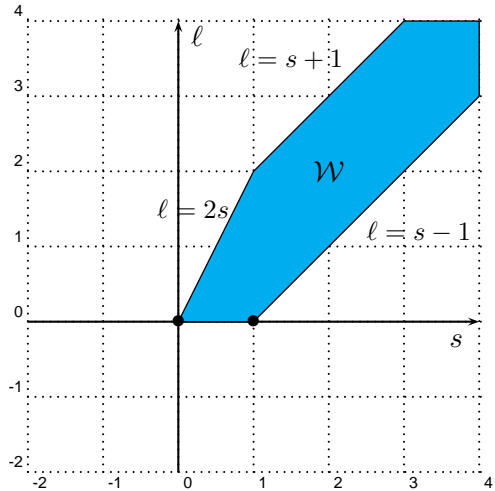

In order to prove our global result and overcome the difficulty caused by the absence of conservation of the  energy of (\ref{SD}), we use a careful control of its derivative  (\ref{Energy-1}) for solutions in $H^1(\R^2)\times L^2(\R^2)$. This method requires the availability of a local theory in this space, a case which is not covered in the previous literature and had to be derived here as well. More precisely, we prove local well-posedness in dimensions $N=2,3$, for initial data in $H^s\times H^{\ell}$ with $s$ and $\ell$ satisfying  $\max \{0, s-1\} \le \ell \le \min\{2s, s+1\}$ (see Theorem \ref{local-theorem} and Figure \ref{Figure-LWP}).

\section{\textbf{Local well-posedness in dimensions $\boldsymbol{N=2, 3}$}}\label{S-LWP}

In this section we obtain local well-posedness for the system (\ref{SD}) following the approach used by Ginibre, Tsutsumi and Velo in \cite{Ginibre}
for the Zakharov system.  As in \cite{Ginibre} we measure the solutions in appropriate Bourgain space, or Fourier restriction, norms. More specifically, we consider
the solution $(u,v)$ of the system (\ref{SD}) in the space $\mathcal{Y}_{s, \ell}=X^{s,\,b}\times  H^{\ell,\, c}$,
the completion of the product of Schwartz spaces  $\mathcal{S}(R^{N+1})\times \mathcal{S}(\R^{N+1})$ with respect to the norm
$$\|(u,v)\|_{\mathcal{Y}_{s, \ell}}=\|u\|_{X^{s,\,b}} + \|v\|_{H^{\ell,\, c}},$$
where
\begin{equation}\label{LWP-norm-v}
\|v\|_{H^{\ell, c}}:=\left \|(1+|\xi|)^{\ell}(1+|\tau|)^c \widehat{v}(\xi,\tau) \right \|_{L^2_{\xi, \tau}},
\end{equation}
and
\begin{equation}\label{LWP-norm-u}
\|u\|_{X^{s,b}}:=
\left\|(1+|\xi|)^s(1+|\tau+ \tfrac{1}{2}|\xi|^2|)^b \widehat{u}(\xi,\tau)\right\|_{L^2_{\xi,\tau}},
\end{equation}
is the Fourier restriction norm associated to the Schr\"odinger group $S(t)$ defined in (\ref{Schrodinger-Group}).
In these definitions $\widehat{f}(\xi,\tau)$ denotes the space-time Fourier transform of $f(x,t)$
and  $|\xi|$  is the Euclidean norm of the frequency vector  $\xi\in \R^N$.

We recall that $X^{s, b} \hookrightarrow C\left(\R; H^s(\R^N)\right)$ and
$ H^{\ell, c} \hookrightarrow C\left (\R; H^{\ell}(\R^N)\right)$ for
all $s, \ell \in \R$ if $b,\, c> 1/2$.

Now we state the main local well-posedness result.

\begin{theorem}\label{local-theorem}
Let $N=2,3$. For any $(u_0,v_0)\in H^s(\R^N) \times H^{\ell}(\R^N)$, with $s$ and $\ell$ satisfying the
conditions:
\begin{equation}\label{local-theorem-a}
\max \{0, s-1\} \le \ell \le \min\{2s, s+1\}
\end{equation}
there exists a positive time $T=T(\|u_0\|_{H^s}, \|v_0\|_{H^{\ell}})$ and
a unique solution $(u(t),v(t))$ of the initial value problem
(\ref{SD}) on the time interval $[0,T]$, such that
\begin{enumerate}
\item [(i)\;] $\left(\psi_T\, u, \psi_T\, v\right)\in X^{s,b}\times H^{\ell, c}$;
\vspace{0.2cm}
\item [(ii)] $(u, v)\in C\left([0,T]; H^s(\R^N) \times H^{\ell}(\R^N)\right)$
\end{enumerate}
for suitable $b$ and $c$ close to $\frac{1}{2}+$ ($\psi_T$ denotes, as usual, a cutoff function for the time interval $[0,T]$).
Moreover, the map $(u_0,v_0) \longmapsto (u(t),v(t))$ is locally
Lipschitz from $H^s(\R^N) \times H^{\ell}(\R^N)$ into $C\left ([0,T];\, H^s(\R^N)
\times H^{\ell}(\R^N)\right)$.
\end{theorem}

\subsection{Preliminary Estimates} In the sequel, we use the following notation.
For $\lambda \in \R$
$$[\lambda]_+=
\begin{cases}
\lambda &\text{if}\quad  \lambda > 0,\\
\varepsilon\; (0< \varepsilon \ll 1) &\text{if}\quad  \lambda = 0,\\
0 & \text{if}\quad  \lambda < 0.
\end{cases}
$$
and we denote by $\lambda \pm$ a number slightly larger, respectively
smaller, than $\lambda $. The bracket $\lan \cdot  \ran$ is defined as $\lan \cdot  \ran = 1+ |\cdot|.$

We introduce the variables
\begin{equation}
\sigma_i=\tau_i + \tfrac{1}{2}|\xi_i|^2,\; \xi_i\in \R^N,\; \tau_i\in \R\; (i=1,2)\quad  \text{and}\quad \sigma=\tau\in \R
\end{equation}
with the convolution structure
\begin{equation}\label{Convolution-Structure}
\xi=\xi_1-\xi_2\quad \text{and}\quad \tau =\tau_1-\tau_2.
\end{equation}
In terms of these variables, the resonance relation for system (\ref{SD})
is the following:
\begin{equation}\label{Dispersive-Relation}
\sigma_1-\sigma_2-\sigma= \tfrac{1}{2}(|\xi_1|^2-|\xi_2|^2).
\end{equation}

\begin{lemma}\label{p.ginibre-s} Let $N=2,3$,\, $b_0=\frac{1}{2}+$,\, $0\le \gamma \le 1$\, and\, $a$,  $a_1$ and $a_2$ be non-negative  numbers
satisfying the conditions
\begin{align}
&(1-\gamma)\max\{a,a_1,a_2\}\le b_0 \le  (1-\gamma)(a+a_1+a_2),\label{Fundamental-Lemma-A}\\
&(1-\gamma)a<b_0.\label{Fundamental-Lemma-B}
\end{align}
Let $m$ be such that
\begin{equation}\label{Fundamental-Lemma-C}
m\geq N/2+1-(1-\gamma)\frac{(a+a_1+a_2)}{b_0}\geq 0
\end{equation}
with strict inequality on the left of (\ref{Fundamental-Lemma-C}) if equality holds on the right of
(\ref{Fundamental-Lemma-A}) or if $a_1=0$.
In addition, let $a'\geq \gamma a,\, a_1'\geq \gamma a_1,\, a_2'\geq \gamma a_2$ and let
$h, h_1, h_2 \in L^2(\R^{N+1})$ be such  that $\cF^{-1}(\langle\sigma\rangle^{-a'}\widehat{h}\,)$,
$\cF^{-1}(\langle\sigma_i\rangle^{-a'_i}\widehat{h_i}\,)\, (i=1,2)$ have support in $|t| \le CT$.
Then, for
\begin{equation}
\theta=\gamma\sum\limits_{j=a,a_1,a_2} j(1-[j'-1/2]_+/j')
\end{equation}
the inequalities
\begin{equation}\label{Fundamental-Lemma-D}
\int \frac{|\widehat{h}(\xi,\tau)\widehat{h}_1(\xi_1,\tau_1) \widehat{h}_2(\xi_2,\tau_2)|}{\langle\sigma\rangle^{a}
\langle\sigma_1\rangle^{a_1}\langle \sigma_2\rangle^{a_2}
\langle\xi\rangle^{m}}\lesssim T^{\theta}\|h\|_{L^2} \|h_1\|_{L^2} \|h_2\|_{L^2}
\end{equation}
and
\begin{equation}\label{Fundamental-Lemma-E}
\int \frac{|\widehat{h}(\xi,\tau)\widehat{h}_1(\xi_1,\tau_1) \widehat{h}_2(\xi_2,\tau_2)|}{\langle\sigma\rangle^{a}
\langle\sigma_1\rangle^{a_1}\langle\sigma_2\rangle^{a_2}
\langle\xi_2\rangle^{m}}\lesssim T^{\theta}\|h\|_{L^2} \|h_1\|_{L^2} \|h_2\|_{L^2},
\end{equation}
hold.
\end{lemma}
\begin{proof}
These estimates follow from Lemma 3.2 in \cite{Ginibre}, for the case
$\sigma_i=\tau_i + |\xi_i|^2$ and $\sigma = \tau \pm|\xi|$, changing  the terms
$\sigma_i\;(i=1,2)$ and  $\sigma$ by $\sigma_i=\tau_i + \tfrac{1}{2}|\xi_i|^2\;(i=1,2)$
and $\sigma=\tau$, respectively.
\end{proof}

\subsection{Bilinear Estimates} It is well known that in the framework of Bourgain spaces local well-posedness results can usually be reduced to the proof of adequate $k$-linear estimates. In the present case, to prove Theorem
\ref{local-theorem} it suffices to establish the following two bilinear estimates:

\begin{proposition}\label{p.uv} Let $s\ge 0$, $\ell \geq \max\{0, s-1\}$ and the functions $u$ and $v$ be supported in time in the region $|t|\le CT$.
Then, the  bilinear estimate
$$\|uv\|_{X^{s,-b_1}}\lesssim T^{\theta}\|u\|_{X^{s,b_2}} \|v\|_{H^{\ell,c}}$$
holds provided $c=\tfrac{1}{2}+\varepsilon$, $b_1=\tfrac{1}{2}-\varepsilon_1$ and $b_2=\tfrac{1}{2}+\varepsilon_2$
for an adequate selection of the parameters $0\le \varepsilon, \varepsilon_1, \varepsilon_2  \ll 1$.
\end{proposition}

\begin{proposition}\label{p.u2} Let $s\ge 0$, $\ell \le \min\{2s, s+1\}$ and the functions $u$ and $w$ be supported in time in the region $|t|\le CT$.
Then, the  bilinear estimate
$$\|(u\cdot\overline{w})\|_{H^{\ell,-b}}\lesssim T^{\theta} \|w\|_{X^{s,b_3}}\|u\|_{X^{s, b_3}}$$
holds provided $b=\tfrac{1}{2}-\varepsilon$, $b_3=\tfrac{1}{2}+\varepsilon_3$
for an adequate selection of the parameters $0\le \varepsilon, \varepsilon_3 \ll 1$.
\end{proposition}

The proofs  of the Propositions \ref{p.uv} and \ref{p.u2} follow similar arguments as the ones used
in \cite{Ginibre} to prove Lemmas 3.4 and 3.5, for the Zakharov system in all dimensions. Thus, we only present here
the proof of Proposition \ref{p.uv}, corresponding to Lemma 3.4  of \cite{Ginibre} in our context (dimensions $N=2,3$), followed by a brief sketch of the proof of Proposition \ref{p.u2}.

\subsection{Proof of the Proposition \ref{p.uv}}
We define
$$
\widehat{h}(\xi, \tau)=\langle \xi \rangle^{\ell}\langle \sigma \rangle^c\widehat{v}(\xi, \tau)\quad
\text{and}
\quad \widehat{h}_2(\xi_2, \tau_2)=\langle \xi_2 \rangle^{s}\langle \sigma_2 \rangle^{b_2}\widehat{u}(\xi_2, \tau_2).
$$
In order to estimate $\|uv\|_{X^{s,-b_1}}$ by duality arguments,  we take the the scalar product
with a generic function in $X^{-s,b_1}$  with Fourier transform
$\langle \xi_1 \rangle^{s}\langle \sigma_1 \rangle^{-b_1}\widehat{h_1}(\xi_1, \tau_1)$ and
$h_1\in L^2({\R^{N+1}})$. Then, the bilinear estimate in
Proposition \ref{p.uv}  takes the form
\begin{equation}\label{S-Estiamte}
|S(h, h_1,h_2)|\lesssim T^{\theta}\|h\|_{L^2}\|h_1\|_{L^2}\|h_2\|_{L^2},
\end{equation}
where
\begin{equation}\label{S-Functional}
S(h, h_1,h_2)=\int\frac{\widehat{h}\widehat{h}_1\widehat{h}_2\langle \xi_1 \rangle^{s}}
{\langle \sigma \rangle^c\langle \sigma_1 \rangle^{b_1}\langle \sigma_2 \rangle^{b_2}\langle \xi_2 \rangle^{s}\langle \xi \rangle^{\ell}}.
\end{equation}
First, we note that if $0\le s \le \ell $ then we have
$$\frac{\langle \xi_1 \rangle^{s}}{\langle \xi_2\rangle^{s}\langle \xi \rangle^{\ell}}\lesssim 1.$$
Then, taking  $(\varepsilon, \varepsilon_1, \varepsilon_2)=(0, \varepsilon, \varepsilon)$ and applying
Lemma \ref{p.ginibre-s}-(\ref{Fundamental-Lemma-D}) with
\begin{align}
&(a',a_1',a_2')=(a,a_1,a_2)=(c, b_1,b_2)=(\tfrac{1}{2},\, \tfrac{1}{2}-\varepsilon,\, \tfrac{1}{2}+\varepsilon),\label{S-First-Estimate-A}\\
&(1-\gamma)=\tfrac{N+2}{3}b_0,\label{S-First-Estimate-B}\\
\intertext{and}
&m=0 \ge N/2+1-(1-\gamma)\frac{(c+b_1+b_2)}{b_0}=0,\label{S-First-Estimate-C}
\end{align}
we obtain
\begin{equation}\label{S-First-Estimate}
|S|\lesssim \int \frac{|\widehat{h}\widehat{h}_1\widehat{h}_2|}
{\langle \sigma \rangle^c\langle \sigma_1 \rangle^{b_1}\langle \sigma_2 \rangle^{b_2}}
\lesssim T^{\theta}\|h\|_{L^2}\|h_1\|_{L^2}\|h_2\|_{L^2},
\end{equation}

To estimate the functional $S$ in the case $s\ge \ell $ we divide the analysis into two cases by considering two
integration  subregions:

\smallskip
\noindent{\textbf{Case 1:} $\boldsymbol{ |\xi_1| \le 2 |\xi_2|}$}. Here, $s\ge 0$ implies
$\langle \xi_1 \rangle^{s}\lesssim \langle \xi_2 \rangle^{s}$; so, the contribution $S_1$
of this subregion to $S$ is given by
\begin{equation}\label{S-Functional-Case1-A}
|S_1|\lesssim \int \frac{|\widehat{h}\widehat{h}_1\widehat{h}_2|}
{\langle \sigma \rangle^c\langle \sigma_1 \rangle^{b_1}\langle \sigma_2 \rangle^{b_2}\langle \xi \rangle^{\ell}}
\lesssim T^{\theta}\|h\|_{L^2}\|h_1\|_{L^2}\|h_2\|_{L^2},
\end{equation}
where we have used Lemma \ref{p.ginibre-s}-(\ref{Fundamental-Lemma-D}) with conditions (\ref{S-First-Estimate-A}),
(\ref{S-First-Estimate-B}) and
\begin{equation}\label{S-Functional-Case1-B}
m=\ell \ge N/2+1-(1-\gamma)\frac{(c+b_1+b_2)}{b_0}=0.
\end{equation}

\smallskip
\noindent{\textbf{Case 2:} $\boldsymbol{ |\xi_1| \ge 2 |\xi_2|}$}. In this situation we have that
$|\xi|\sim |\xi_1|$. Also, from the resonance relation (\ref{Dispersive-Relation})
it follows that
\begin{equation}\label{S-Functional-Case2-A}
|\sigma_1-\sigma_2-\sigma|= \tfrac{1}{2}\Bigl||\xi_1|^2-|\xi_2|^2\Bigl | \ge \tfrac{3}{8}|\xi_1|^2
\Longrightarrow |\xi_1|^2 \le \tfrac{9}{8}\max\{|\sigma|, |\sigma_1|, |\sigma_2|\}.
\end{equation}
Then, in view of (\ref{S-Functional-Case2-A}) and using the fact that $s\ge \ell$ we estimate the contribution $S_2$ of this subregion to $S$ by
\begin{equation}\label{S-Functional-Case2-B}
|S_2|\lesssim \int\frac{|\widehat{h}\widehat{h}_1\widehat{h}_2|\langle \xi_1 \rangle^{s-\ell}}
{\langle \sigma \rangle^c\langle \sigma_1 \rangle^{b_1}\langle \sigma_2 \rangle^{b_2}\langle \xi_2 \rangle^{s}}
\lesssim \int\frac{|\widehat{h}\widehat{h}_1\widehat{h}_2|\langle \sigma^* \rangle^{\rho(s,\ell)}}
{\langle \sigma \rangle^c\langle \sigma_1 \rangle^{b_1}\langle \sigma_2 \rangle^{b_2}\langle \xi_2 \rangle^{s}},
\end{equation}
where $\rho(s,\ell)=\frac{s-\ell}{2}$ and $\sigma^* = \max \bigl\{\langle \sigma \rangle, \langle \sigma_1 \rangle, \langle \sigma_2 \rangle \bigl\}$.

Now we consider the condition
\begin{equation}\label{S-Functional-Case2-c}
\rho(s,\ell) =\frac{s-\ell}{2} \le \min\{c, b_1, b_2\}=\frac{1}{2}-\varepsilon_1 \Longleftrightarrow s-1 + 2\varepsilon_1 \le \ell ,
\end{equation}
which guarantees that $b-\rho$, $b_1-\rho$ and $b_2-\rho$ are nonnegative.
Next,  we establish  conditions that  allow us to apply Lemma \ref{p.ginibre-s}-(\ref{Fundamental-Lemma-D}), with
\begin{equation}\label{S-Functional-Case2-C}
(a',a_1',a_2')=(a, a_1, a_2)=
\begin{cases}
(c-\rho,\, b_1,\, b_2)& \text{if}\quad \sigma^*=\sigma,\\
(b,\, b_1-\rho,\, b_2)& \text{if}\quad \sigma^*=\sigma_1,\\
(b,\, b_1,\, b_2-\rho)& \text{if}\quad \sigma^*=\sigma_2,\\
\end{cases}
\end{equation}
and $m=s$, to obtain the desired estimate
\begin{equation}
|S_2|\lesssim \int\frac{|\widehat{h}\widehat{h}_1\widehat{h}_2|\langle \sigma^* \rangle^{\rho}}
{\langle \sigma \rangle^c\langle \sigma_1 \rangle^{b_1}\langle \sigma_2 \rangle^{b_2}\langle \xi_2 \rangle^{s}}
\lesssim T^{\theta}\|h\|_{L^2}\|h_1\|_{L^2}\|h_2\|_{L^2}.
\end{equation}
For this purpose, it suffices to take $0< \gamma <1$ such that
\begin{align}
&b_0\le (1-\gamma)(c+b_1+b_2-\rho),\label{S-Functional-Case2-D}\\
&s\ge  N/2+1-(1-\gamma)\frac{(c+b_1+b_2-\rho)}{b_0}\ge 0.\label{S-Functional-Case2-E}
\end{align}
Choosing $\gamma'=\gamma\frac{c+b_1+b_2-\rho}{c+b_1+b_2}\in (0,1)$ conditions (\ref{S-Functional-Case2-D})-(\ref{S-Functional-Case2-E})
are equivalent to taking
$\gamma'$ such that
\begin{align}
&b_0\le (1-\gamma')(c+b_1+b_2)-\rho(s,\ell),\label{S-Functional-Case2-F}\\
&s\ge  N/2+1-(1-\gamma')\frac{(c+b_1+b_2)}{b_0}+\frac{\rho(s,\ell)}{b_0}\ge 0.\label{S-Functional-Case2-G}
\end{align}
Now, we take  $(\varepsilon, \varepsilon_1, \varepsilon_2)=(0, \varepsilon, \varepsilon)$ and $\gamma'$ satisfying $(1-\gamma')=\tfrac{N+2}{3}b_0$
and then from (\ref{S-Functional-Case1-B})
we have
$$ s >  \ell +\frac{\rho(s,\ell)}{b_0}\ge \underbrace{N/2+1-(1-\gamma')\frac{(c+b_1+b_2)}{b_0}}_{=0}+\frac{\rho(s,\ell)}{b_0}\ge 0$$
as desired. This completes the proof.

\subsection{Proof of the Proposition \ref{p.u2}} Following the same ideas as in the proof of Proposition \ref{p.uv}, in this case
we need to estimate the following functional:
\begin{equation}\label{W-Functional}
W(h, h_1,h_2)=\int\frac{\widehat{h}\widehat{h}_1\widehat{h}_2\langle \xi \rangle^{\ell}}
{\langle \sigma \rangle^b\langle \sigma_1 \rangle^{b_3}\langle \sigma_2 \rangle^{b_3}\langle \xi_1 \rangle^{s}\langle \xi_2 \rangle^{s}}.
\end{equation}
Unlike the Zakharov system, here we do not have the presence of the derivative term ($|\xi|$)
in the numerator of $W$. Thus,  we estimate (\ref{W-Functional}) in the same way as in the proof of
Lemma 3.5 in \cite{Ginibre} ($N=2, 3$), replacing $|\xi|\langle \xi \rangle^{\ell}$  by $\langle \xi \rangle^{\ell}$,
which is equivalent to changing $\ell +1$ by $\ell$ in the internal computations.

\begin{remark} In the case $|\ell-s|<1$ the bilinear estimates in Propositions \ref{p.uv} and \ref{p.u2}
hold for small positive numbers  $\varepsilon, \varepsilon_1, \varepsilon_2$ and $\varepsilon_3$ which allows
taking $b=\frac{1}{2}+$ and $c=\frac{1}{2}+$ in Theorem \ref{local-theorem}, so that the
corresponding immersions $X^{s, b} \hookrightarrow C\left(\R; H^s(\R^N)\right)$ and
$ H^{\ell, b}\hookrightarrow C\left (\R; H^{\ell}(\R^N)\right)$ are guaranteed. On the other hand,
if $|\ell-s|=1$, one must  take $b=c=1/2$ in Theorem \ref{local-theorem}; then, to guarantee  the required
immersions, we need to establish extra bilinear estimates in the norms
$$\|u\|_{\tilde{X}^s}:=\left \| \langle \xi \rangle ^s \langle \tau+ \tfrac{1}{2}|\xi|^2\rangle^{-1} \widehat{u}(\xi,\tau)\right\|_{L^2_{\xi}L^1_{\tau}}
$$
and
$$\|v\|_{\tilde{H}^{\ell}}:=\left \| \langle \xi \rangle^{\ell} \langle \tau \rangle^{-1}
\widehat{v}(\xi,\tau) \right\|_{L^2_{\xi}L^1_{\tau}},$$
which follow in the same manner as in the proof of Lemmas 3.6 and 3.7 in \cite{Ginibre}.
\end{remark}

\subsection{Proof of the Theorem \ref{local-theorem}} The proof follows the, now standard, contraction method applied
to a localized in time cut-off integral formulation associated to the  system (\ref{SD}), in the Bourgain spaces defined by the norms \eqref{LWP-norm-v} and \eqref{LWP-norm-u} (see \cite{Ginibre}, for example, for complete details of a similar proof). As is well known, the success of this method relies almost exclusively on the availabilty of certain multilinear estimates in these norms for the nonlinear terms of the equations. In our case these estimates are the ones obtained in Propositions \ref{p.uv} and \ref{p.u2}.

We start with the following integral system  \footnote{Observe that for the $v$ equation \eqref{int-cut-2}, we are not 
using the standard Duhamel, or variation of parameters, form as in \eqref{eq. de v}. Both forms could be used here, but our
choice makes computations slightly simpler.} 
\begin{align}
u(\cdot, t)&=\Phi_1(u,v):=\psi_1(t)S(t)u_0-i\psi_T(t)\int_0^tS(t-t')\psi_{2T}^2(t')u(\cdot, t')v(\cdot, t')dt', \label{int-cut-1}\\
v(\cdot, t)&=\Phi_2(u,v):=\psi_1(t)v_0+ \psi_T(t)\int_0^t\Bigl [\tfrac{\lambda}{\mu}|\psi_{2T}(t')u(\cdot, t')|^2 -\psi_{2T}(t')v(\cdot, t')\Bigl]dt'.\label{int-cut-2}
\end{align}
Here $S(t)$ is given by \eqref{Schrodinger-Group}, $\psi_1 \in C^{\infty}(\R;\; \R^+)$ is an even function, such that $0\leq \psi_1\leq 1$ and
\begin{equation}\label{cut-1}
\psi_1(t) =
\begin{cases}
1, \quad |t|\leq 1,\\
0, \quad |t|\geq 2
\end{cases}
\end{equation}
and  $\psi_T(t) = \psi_1(t/T)$ for $0\leq T\leq 1$.

Now, choosing $T$ sufficiently small,  we solve (\ref{int-cut-1})-(\ref{int-cut-2}) by contraction  in
the complete metric space
$$ \Sigma= \Bigl\{(u,v)\in X^{s,b}\times H^{\ell, c}; \; \|u\|_{X^{s,b}}\le N_1\; \text{and}\; \|v\|_{H^{\ell, c}} \le N_2 \Bigl\},$$
with $b, c =\frac{1}{2}+$ and the induced norm $\|(u,v)\|_{\Sigma}:= \|u\|_{X^{s,b}} + \|v\|_{H^{\ell, c}},$
where  the positive constants $N_j, j=1,2,$ will be chosen below.

Using the properties of spaces $X^{s,b}$ and $H^{\ell, c}$ (see  (2.19) and Lemma 2.1  in \cite{Ginibre}) and
the bilinear estimates obtained in Propositions \ref{p.uv} and \ref{p.u2} we get
\begin{equation}\label{PLT-A}\begin{split}
\|\Phi_1(u,v)\|_{X^{s,b}}& \le c_0\|u_0\|_{H^s} + c_1T^{1-b-a_1}\|\psi_{2T}^2uv\|_{X^{s,-a_1}},\\
&\le c_0\|u_0\|_{H^s} + c_1T^{\theta_1}\|u\|_{X^{s,b}} \|v\|_{H^{\ell,c}}\\
&\le c_0\|u_0\|_{H^s} + c_1T^{\theta_1}N_1N_2
\end{split}\end{equation}
and
\begin{equation}\label{PLT-B}\begin{split}
\|\Phi_2(u,v)\|_{H^{\ell,c}}& \le c_0\|v_0\|_{H^{\ell}} + c_2T^{1-c-a_2}\left ( \|\psi_{2T}^2|u|^2\|_{H^{\ell,-a_2}} \|\psi_{2T}^2v\|_{H^{\ell,-a_2}}\right)\\
&\le c_0\|v_0\|_{H^{\ell}} + c_2T^{\theta_2}\left (\|u\|^2_{X^{s, b}}+  \|v\|_{H^{\ell,c}} \right )\\
&\le c_0\|v_0\|_{H^{\ell}} + c_2T^{\theta_2}(N_1^2+N_2),
\end{split}\end{equation}
where $-\frac{1}{2}< -a_1\le 0\le b\le 1-a_1$,\; $-\frac{1}{2}< -a_2\le 0\le c\le 1-a_2$ and $\theta_i,\; i=1,2$ are positive.

Taking $N_1=2c_0\|u_0\|_{H^s}$ and $N_2=2c_0\|v_0\|_{H^{\ell}}$ from (\ref{PLT-A}) and (\ref{PLT-B})
it follows that $\Phi(u,v)=(\Phi_1(u,v),\,\Phi_2(u,v))\in \Sigma$ for $(u,v)\in \Sigma$ and for small enough $T$ satisfying
$$T^{\theta_1}\le \frac{1}{2c_1N_2}\quad \text{and}\quad  T^{\theta_2}\le \frac{N_2}{2c_2(N_1^2+N_2)}.$$

The contraction condition can be obtained in a similar way and the proof is finished.

\begin{remark} In the case $b=c=1/2$ the proof of the Theorem \ref{local-theorem} follows by similar arguments, but using
the following modified norms
$$\|\!|u\|\!|_s:=\|u\|_{X^{s,1/2}}+ \left \| \langle \xi \rangle ^s \langle \tau+ \tfrac{1}{2}|\xi|^2\rangle^{-1} \widehat{u}(\xi,\tau)\right\|_{L^2_{\xi}L^1_{\tau}}
$$
and
$$\|\!|v\|\!|_{\ell}:=\|u\|_{H^{\ell,1/2}}+\left \| \langle \xi \rangle^{\ell} \langle \tau \rangle^{-1} \widehat{v}(\xi,\tau) \right\|_{L^2_{\xi}L^1_{\tau}}$$
in order to obtain the immersions in the spaces $C\left(\R; H^s(\R^N)\right)$ and $C\left (\R; H^{\ell}(\R^N)\right)$.
\end{remark}

\section{\textbf{Global well-posedness for the critical model}}\label{S-APE}

As mentioned in the introduction, in dimension $N=2$, the system (\ref{SD}) is a perturbation of the  scaling-critical cubic NLS equation (\ref{NLS}).
In this section we derive {\it a priori} estimates in the energy space $H^1(\R^2)\times L^2(\R^2)$ for the focusing and defocusing cases
of (\ref{SD}), which allow us to extend the local solutions obtained in the previous section to all positive times.

\begin{theorem} \label{global-theorem}
Let $(u_0,v_0)\in H^1(\R^2)\times L^2(\R^2)$.
Then, for all $T>0$, there exists a unique solution
$$(u,v)\in C\left([0,T];\; H^{1}(\R^2)\times L^{2}(\R^2)\right)$$
to the Initial Value Problem (\ref{SD}).
\end{theorem}

\begin{proof}
In view of the local well-posedness result detailed in the last section and the conservation of
the $L^2$-norm of $u$, we only need to obtain an {\it a priori} bound for the function
\[f(t):= \|\nabla u(\cdot, t)\|_{L^2}^2 +\|v(\cdot,t)\|_{L^2}^2.\]

We begin by estimating $\|v(\cdot,t)\|_{L^2}^2.$  Using the explicit representation for $v$, given by
\begin{equation}\label{APE-1}
v(x,t)= e^{-t/\mu}v_0(x) + \tfrac{\lambda}{\mu}\int_0^t\,
e^{-(t-t')/\mu}|u(x,t')|^2\,dt',
\end{equation}
and applying the Minkowski and Gagliardo-Nirenberg inequalities we get
\begin{equation}\begin{split}\label{APE-2}
\|v(\cdot, t)\|_{L^2}&\le \|v_0\|_{L^2} + \tfrac{1}{\mu}
\int_0^t\,e^{-(t-t')/\mu}\|u(\cdot,t')\|_{L^4}^2\,dt'\\
&\le \|v_0\|_{L^2} + \tfrac{\beta^2}{\mu}
\int_0^t\,\|u(\cdot, t')\|_{L^2}\|\nabla u(\cdot, t')\|_{L^2}\,dt'\\
&=\|v_0\|_{L^2} + \tfrac{\beta^2\|u_0\|_{L^2}}{\mu}\int_0^t\,\|\nabla u(\cdot, t')\|_{L^2}\,dt',
\end{split}\end{equation}
where $\beta$ is the constant from the Gagliardo-Nirenberg inequality.
Now, we use H\"older's inequality to  obtain, from (\ref{APE-2})
\begin{equation}\begin{split}\label{APE-3}
\|v(\cdot, t)\|^2_{L^2} &\le 2\|v_0\|^2_{L^2} + 2\left(\tfrac{\beta^2\|u_0\|_{L^2}}{\mu}\int_0^t\,\|\nabla u(\cdot, t')\|_{L^2}\,dt'\right)^2\\
&\le 2\|v_0\|^2_{L^2} + \tfrac{2\beta^4\|u_0\|^2_{L^2}}{\mu^2}\,t\,\int_0^t\,\|\nabla u(\cdot, t')\|^2_{L^2}\,dt'\\
&\le 2\|v_0\|^2_{L^2} + \tfrac{2\beta^4\|u_0\|^2_{L^2}}{\mu^2}\,t\,\int_0^t\,f(t')\,dt'.
\end{split}\end{equation}

On the other hand, from (\ref{Energy-2}), we have
\begin{equation}\begin{split}\label{APE-4}
\|\nabla u(\cdot, t)\|^2_{L^2}&= E(t) -2\int_{\R^2}v(\cdot, t)|u(\cdot, t)|^2dx+\lambda\|v(\cdot, t)\|^2_{L^2}\\
&\le|E(t)|+ \left |2\int_{\R^2}v(\cdot, t)|u(\cdot, t)|^2dx\right|+\|v(\cdot, t)\|^2_{L^2}.
\end{split}\end{equation}

We begin by treating the term $\displaystyle \int 2v|u|^2dx$, which, in view of (\ref{APE-1}), can be rewritten in the form
\begin{equation}\label{APE-5}
2\int_{\R^2}v(\cdot, t)|u(\cdot, t)|^2dx = A(t) +B(t),
\end{equation}
where
\[
A(t)=2e^{-t/\mu}\int_{\R^2}v_0|u|^2dx\quad \text{and}\quad B(t)=\tfrac{2\lambda}{\mu}\int_{\R^2}|u(x,t)|^2\int_0^te^{-(t-t')/\mu}|u(x,t')|^2dt'dx.
\]
We now proceed with the estimates of the terms $A(t)$ and $B(t)$.  To estimate
$A(t)$ we use the H\"older and Gagliardo-Nirenberg inequalities to obtain
\begin{equation}\begin{split}\label{APE-6}
|A(t)|&\le 2e^{-t/\mu}\|v_0\|_{L^2}\|u(\cdot,t)\|^2_{L^4}\\
&\le 2\|v_0\|_{L^2}\beta^2\|u(\cdot, t)\|_{L^2}\|\nabla u(\cdot,t)\|_{L^2}\\
&=2 \|v_0\|_{L^2}\beta^2\|u_0\|_{L^2}\|\nabla u(\cdot,t)\|_{L^2}\\
&\le 4\|v_0\|^2_{L^2}\beta^4\|u_0\|^2_{L^2} + \tfrac{1}{4}\|\nabla u(\cdot,t)\|_{L^2}^2.
\end{split}\end{equation}
Now we estimate $B(t)$. By the H\"older, Minkowski and Gagliardo-Nirenberg inequalities we get
\begin{equation}\begin{split}\label{APE-7}
|B(t)|&\le\tfrac{2}{\mu}\|u(\cdot, t)\|_{L^4}^2\int_0^te^{-\frac{(t-t')}{\mu}}\|u(\cdot, t')\|_{L^4}^2dt'\\
&\le\tfrac{1}{\mu}\int_0^t2\|u(\cdot, t)\|_{L^4}^2\|u(\cdot, t')\|_{L^4}^2dt'\\
&\le \tfrac{1}{\mu}\int_0^t\|u(\cdot, t)\|_{L^4}^4dt'+ \tfrac{1}{\mu}\int_0^t\|u(\cdot, t')\|_{L^4}^4dt'\\
&\le \tfrac{\beta^4}{\mu}\,t\,\|u_0\|_{L^2}^2 \|\nabla u(\cdot,t)\|^2_{L^2}+\tfrac{\beta^4}{\mu}\|u_0\|_{L^2}^2\int_0^t
\|\nabla u(\cdot,t')\|_{L^2}^2dt'\\
&\le \tfrac{1}{4}\|\nabla u(\cdot,t)\|^2_{L^2}+ \tfrac{\beta^4}{\mu}\|u_0\|_{L^2}^2\int_0^t f(t')dt',
\end{split}\end{equation}
for all $t\in [0,\,T_{\mu}]$, where  $T_{\mu}=\frac{\mu}{4\beta^4\|u_0\|_{L^2}^2}$. Then, (\ref{APE-6}) and (\ref{APE-7}) yield
\begin{equation}\label{APE-8}
\left |2\int_{\R^2}v(\cdot, t)|u(\cdot, t)|^2dx\right|\le  4\|v_0\|^2_{L^2}\beta^4\|u_0\|^2_{L^2}
+ \tfrac 12 \|\nabla u(\cdot,t)\|_{L^2}^2 + \tfrac{\beta^4}{\mu}\|u_0\|_{L^2}^2\int_0^t f(t')dt',
\end{equation}
for all $t\in [0,\,T_{\mu}]$.

Next, we estimate the growth of $E(t)$. Using \eqref{Energy-1}, we have
\begin{equation}\begin{split}\label{APE-9}
|E(t)|&=\displaystyle \left|E_0+2\lambda \mu\int_0^t\left(\int_{\R^2}v_t^2(x,t')dx\right)dt'\right|\\
&\leq |E_0|+2\mu\int_0^t\left(\int_{\R^2} v_t^2(x,t')dx\right)dt'\\
&\leq |E_0|+\frac 2{\mu}\int_0^t\left(\int_{\R^2} \bigl(\lambda|u(x,t')|^2-v(x,t')\bigl)^2dx\right)dt' \quad \textrm { (by (\ref{SD}))}\\
&\leq |E_0|+\frac 4{\mu}\int_0^t\left(\|u(\cdot, t')\|^4_{L^4}+\|v(\cdot, t')\|_{L^2}^2\right)dt'\\
&\leq |E_0|+\frac 4{\mu}\int_0^t\left(\beta^4\|u_0\|^2_{L^2}\|\nabla u(t')\|^2_{L^2}+\|v(t')\|_{L^2}^2\right)dt'\\
&\leq |E_0|+\frac 4{\mu}\left(\beta^4\|u_0\|_{L^2}^2+1\right)\int_0^tf(t')dt'.
\end{split}\end{equation}
Collecting the information in (\ref{APE-4}), (\ref{APE-8}) and (\ref{APE-9}) we have
\begin{equation*}\begin{split}\label{APE-10}
\|\nabla u(\cdot,t)\|_{L^2}^2&\le 2|E_0|+ 8\|v_0\|^2_{L^2}\beta^4\|u_0\|^2_{L^2} + 2\|v(\cdot, t)\|^2_{L^2}\\
&\hspace{2.5cm} +\;\tfrac {2}{\mu}\left(5\beta^4\|u_0\|_{L^2}^2+4\right)\int_0^tf(t')dt',
\end{split}\end{equation*}
for all $t\in \bigl[0,\, T_{\mu}\bigl]$, from which follows
\begin{equation}\begin{split}\label{APE-11}
f(t)&\le 2|E_0|+ 8\|v_0\|^2_{L^2}\beta^4\|u_0\|^2_{L^2} + 3\|v(\cdot, t)\|^2_{L^2}\\
&\hspace{2.5cm} +\;\tfrac {2}{\mu}\left(5\beta^4\|u_0\|_{L^2}^2+4\right)\int_0^tf(t')dt',
\end{split}\end{equation}
for all $t\in \bigl[0,\, T_{\mu}\bigl]$.

Finally, we combine (\ref{APE-3}) with (\ref{APE-11}) to obtain
\begin{equation}\label{APE-12}
f(t)\leq \alpha_0 + \alpha_1\int_0^tf(t')dt', \quad \text{for all}\quad t\in \bigl[0,\, T_{\mu}\bigl],
\end{equation}
with
\begin{align}
&\alpha_0 = 2|E_0|+ 4\|v_0\|^2_{L^2}\left( 2\beta^4\|u_0\|^2_{L^2} +\tfrac 32\right),\nonumber\\
&\alpha_1=\tfrac {2}{\mu}\left(5\beta^4\|u_0\|_{L^2}^2+\tfrac{19}4\right)\nonumber.
\end{align}
Hence, by Gronwall's lemma,
$$f(t)\le \alpha_0e^{\alpha_1 t},\quad t\in \bigl[0,\, T_{\mu}\bigl].$$

 Since the time $T_{\mu}=\frac{\mu}{4\beta^4\|u_0\|_{L^2}^2}$ depends only on the conserved quantity $\|u_0\|_{L^2}$, we can iterate this procedure in order to extend this solution to all positive times. Note, however, that the solution can blow-up at infinity.
\end{proof}
\section{Concluding Remarks}

\subsection{Global well-posedness in $\boldsymbol{L^2\times L^2\, (N=2,3)}$}  We observe that following the same
ideas outlined in Remark 5.5 in \cite{Corcho-Matheus} we can extend our local results in $L^2\times L^2$
to any positive time $T$.

\subsection{Global well-posedness in $\boldsymbol{H^1\times H^1\, (N=1)}$}
The results in \cite{Corcho-Matheus} concerning local well-posedness in one dimension do not include the case $(u_0,v_0)\in H^1\times L^2$. However, they include $(u_0,v_0)\in H^1\times H^1$. Our proof in Section 3 can be adapted in order to obtain global well-posedness in this situation, as well.
Indeed, putting
$$g(t)=f(t)+\|v_x(\cdot, t)\|_{L^2}=\|v(\cdot, t)\|_{L^2}^2+\|u_x(\cdot, t)\|_{L^2}^2+\|v_x(\cdot, t)\|_{L^2},$$
and using the following inequality
$$\|u\|_{L^4}^2 \leq \|u\|_{L^{\infty}} \|u\|_{L^2} \leq \|u\|_{H^1} \|u\|_{L^2},$$
instead of the one dimensional Gagliardo-Nirenberg,
we obtain as in (\ref{APE-12})
$$f(t)\leq \tilde{\alpha}_0 + \tilde{\alpha}_1\int_0^tf(t')dt', \quad \text{for all}\quad t\in \bigl[0,\, \tilde{T}\bigl],$$
for some $\tilde{T}=\tilde{T}(\|u_0\|_{L^2},\mu,\beta)$ and where $\tilde{\alpha}_i$, $i=0,1$, depend exclusively on the initial data.

Furthermore, differentiating (\ref{APE-1}) with respect to $x$ and taking the $L^2$-norm yields
\begin{equation}
\begin{split}
\|v_x(\cdot, t)\|_{L^2}&\le \|{v_0}_x\|_{L^2}+ \tfrac{2}{\mu}\int_0^t\|u(t')\|_{L^{\infty}}\|u_x(t')\|_{L^2}dt'\\
&\le \|{v_0}_x\|_{L^2}+\tfrac{2C}{\mu}\int_0^t\|u(t')\|_{H^1}\|u_x(t')\|_{L^2}dt'\\
&\le \|{v_0}_x\|_{L^2}+\tfrac{2C}{\mu}\left(\tfrac{\tilde{T}}{2}\|u_0\|_{L^2}^2+\tfrac 32\int_0^tf(t')dt'\right),
\end{split}
\end{equation}
where we have used the Sobolev embedding $H^1(\R)\hookrightarrow L^{\infty}(\R)$.
Finally, we obtain
$$g(t)\leq \tilde{\alpha}_0+ \|{v_0}_x\|_{L^2}+\tfrac{C\tilde{T}}{\mu}\|u_0\|_{L^2}^2+ \left(\tfrac{3C}{\mu}+\tilde{\alpha}_1\right)\int_0^tg(t')dt', \quad \text{for all}\quad t\in \bigl[0,\, \tilde{T}\bigl].$$
We conclude as in Theorem \ref{global-theorem}.

\subsection{On the blow-up in $\mathbf {H^1\times H^1\, (N=2)}$} Since, in two dimensions, $\|u\|_{H^1}$ does
not control $\|u\|_{L^{\infty}}$, the considerations in the previous Remark do not apply. However, our global result
for initial data $(u_0,v_0)\in H^1\times L^2$ shows that a possible blow-up in $H^1\times H^1$ can only occur for $\|\nabla v\|_{L^2}$.

\subsection{Comparison between the cubic NLS and the Schr\"odinger-Debye equations}
In the next table, we summarize all known results concerning the local and global wellposedness for these equations. It illustrates
the regularization induced by the delay $\mu$ in the Schr\"odinger-Debye system.

\bigskip
\begin{center}
\begin{tabular}{|c|l|l|}
  \hline
  $N$ &
  {\tiny \bf Cubic NLS ($\boldsymbol{H^s}$)} & {\tiny \bf Schr\"odinger-Debye $\boldsymbol{(H^s\times H^{\ell})}$} \\
  \hline
  \hline
  1  &{\tiny l.w.p: $s\ge 0$ (\cite{Ginibrevelo}, \cite{Y-Tsutsumi}, \cite{Cazenave})}& {\tiny l.w.p: $|s|-\frac 12\le \ell <\min\{s+\frac 12; 2s+\frac 12\}$ (\cite{Corcho-Matheus})}\\
  \cline {2-3}
   & {\tiny g.w.p: $s \geq 0$ (Conservation of $L^2$ norm)} & {\tiny g.w.p: $-\frac 3{14}\leq \ell =s \leq 0$ (\cite{Corcho-Matheus}) or $(s,\ell)=(1,1)$} \\
  \hline
  2 &{\tiny l.w.p: $s\ge 0$ (\cite{Ginibrevelo}, \cite{Cazenavecritical}, \cite{Cazenave}) }& {\tiny l.w.p: $\max \{0, s-1\} \le \ell \le \min\{2s, s+1\}$}\\
  \cline {2-3}
   & {\tiny g.w.p ($\lambda=1$): $s > \frac 13$ (\cite{global2})} &  {\tiny g.w.p ($\lambda=\pm 1$): $(s,\ell)=(1,0)$ and $(s,\ell)=(0,0)$} \\
    & {\tiny g.w.p ($\lambda=\pm 1$): $s\ge 0$ for small $L^2$ norm  (\cite{Cazenavecritical})} &\\
  \hline
  3  &{\tiny l.w.p: $s\ge \frac 12$ (\cite{Ginibrevelo}, \cite{Cazenave})}&{\tiny l.w.p: $\max \{0, s-1\} \le \ell \le \min\{2s, s+1\}$}\\
  \cline {2-3}
   & {\tiny g.w.p ($\lambda=1$): $s > \frac 45$ (\cite{global3})} &  {\tiny g.w.p ($\lambda=\pm 1$):  $(s,\ell)=(0,0)$} \\
    & {\tiny g.w.p ($\lambda=\pm 1$): $s\ge \frac 12$ for small $H^{1/2}$ norm  (\cite{Cazenavecritical})} &\\
  \hline
\end{tabular}
\end{center}
\bigskip
{\bf Acknowledgment.} Part of this research was carried while A. J. Corcho was visiting the Department
of Mathematics of IST/Lisbon, Portugal supported by CAPES-Brazil and by the 2010 FCT-CAPES joint project 
{\it Nonlinear waves and dispersion}.

\end{document}